\catcode`\^^Z=9
\catcode`\^^M=10
\expandafter\ifx\csname mdqon\endcsname\relax
\else \endinput \fi
\message{Document Style Option `german'  Version 2 as of 16 May 1988}
%
\ifx\protect\undefined
\let\protect=\relax \fi
{\catcode`\@=11 
\gdef\allowhyphens{\penalty\@M \hskip\z@skip}

\newcount\U@C\newbox\U@B\newdimen\U@D
\gdef\umlauthigh{\def\"##1{{\accent127 ##1}}}
\gdef\umlautlow{\def\"{\protect\newumlaut}}
\gdef\newumlaut#1{\leavevmode\allowhyphens
     \vbox{\baselineskip\z@skip \lineskip.25ex
     \ialign{##\crcr\hidewidth
     \setbox\U@B\hbox{#1}\U@D .01\p@\U@C\U@D
     \U@D\ht\U@B\advance\U@D -1ex\divide\U@D \U@C
     \U@C\U@D\U@D\the\fontdimen1\the\font
     \multiply\U@D \U@C\divide\U@D 100\kern\U@D
     \vbox to .20ex  
     {\hbox{\char127}\vss}\hidewidth\crcr#1\crcr}}\allowhyphens}
\gdef\highumlaut#1{\leavevmode\allowhyphens
     \accent127 #1\allowhyphens}
\gdef\set@low@box#1{\setbox1\hbox{,}\setbox\z@\hbox{#1}\dimen\z@\ht\z@
     \advance\dimen\z@ -\ht1
     \setbox\z@\hbox{\lower\dimen\z@ \box\z@}\ht\z@\ht1 \dp\z@\dp1}
%
\gdef\@glqq{{\ifhmode \edef\@SF{\spacefactor\the\spacefactor}\else
     \let\@SF\empty \fi \leavevmode
     \set@low@box{''}\box\z@\kern-.04em\allowhyphens\@SF\relax}}
\gdef\glqq{\protect\@glqq}
\gdef\@grqq{\ifhmode \edef\@SF{\spacefactor\the\spacefactor}\else
     \let\@SF\empty \fi \kern-.07em``\kern.07em\@SF\relax}
\gdef\grqq{\protect\@grqq}
\gdef\@glq{{\ifhmode \edef\@SF{\spacefactor\the\spacefactor}\else
     \let\@SF\empty \fi \leavevmode
     \set@low@box{'}\box\z@\kern-.04em\allowhyphens\@SF\relax}}
\gdef\glq{\protect\@glq}
\gdef\@grq{\kern-.07em`\kern.07em}
\gdef\grq{\protect\@grq}
\gdef\@flqq{\ifhmode \edef\@SF{\spacefactor\the\spacefactor}\else
     \let\@SF\empty \fi
     \ifmmode \ll \else \leavevmode
     \raise .2ex \hbox{$\scriptscriptstyle \ll $}\fi \@SF\relax}
\gdef\flqq{\protect\@flqq}
\gdef\@frqq{\ifhmode \edef\@SF{\spacefactor\the\spacefactor}\else
     \let\@SF\empty \fi
     \ifmmode \gg \else \leavevmode
     \raise .2ex \hbox{$\scriptscriptstyle \gg $}\fi \@SF\relax}
\gdef\frqq{\protect\@frqq}
\gdef\@flq{\ifhmode \edef\@SF{\spacefactor\the\spacefactor}\else
     \let\@SF\empty \fi
     \ifmmode < \else \leavevmode
     \raise .2ex \hbox{$\scriptscriptstyle < $}\fi \@SF\relax}
\gdef\flq{\protect\@flq}
\gdef\@frq{\ifhmode \edef\@SF{\spacefactor\the\spacefactor}\else
     \let\@SF\empty \fi
     \ifmmode > \else \leavevmode
     \raise .2ex \hbox{$\scriptscriptstyle > $}\fi \@SF\relax}
\gdef\frq{\protect\@frq}

\global\let\original@ss=\ss
\gdef\newss{\leavevmode\allowhyphens\original@ss\allowhyphens}
\global\let\ss=\newss
\global\let\original@three=\3 
%
%
\gdef\german@dospecials{\do\ \do\\\do\{\do\}\do\$\do\&%
  \do\#\do\^\do\^^K\do\_\do\^^A\do\%\do\~\do\"}
\gdef\german@sanitize{\@makeother\ \@makeother\\\@makeother\$\@makeother\&%
\@makeother\#\@makeother\^\@makeother\^^K\@makeother\_\@makeother\^^A%
\@makeother\%\@makeother\~\@makeother\"}
\global\let\original@dospecials\dospecials
\global\let\dospecials\german@dospecials
\global\let\original@sanitize\@sanitize
\global\let\@sanitize\german@sanitize
\gdef\mdqon{\let\dospecials\german@dospecials
        \let\@sanitize\german@sanitize\catcode`\"\active}
\gdef\mdqoff{\catcode`\"12\let\original@dospecials\dospecials
        \let\@sanitize\original@sanitize}

{\mdqoff
\gdef\@UMLAUT{\"}
\gdef\@MATHUMLAUT{\mathaccent"707F}
\gdef\@SS{\mathchar"7019}
\gdef\dq{"}
}
{\mdqon
\gdef"#1{\if\string#1`\glqq{}%
\else\if\string#1'\grqq{}%
\else\if\string#1a\ifmmode\@MATHUMLAUT a\else\@UMLAUT a\fi
\else\if\string#1o\ifmmode\@MATHUMLAUT o\else\@UMLAUT o\fi
\else\if\string#1u\ifmmode\@MATHUMLAUT u\else\@UMLAUT u\fi
\else\if\string#1A\ifmmode\@MATHUMLAUT A\else\@UMLAUT A\fi
\else\if\string#1O\ifmmode\@MATHUMLAUT O\else\@UMLAUT O\fi
\else\if\string#1U\ifmmode\@MATHUMLAUT U\else\@UMLAUT U\fi
\else\if\string#1e\ifmmode\@MATHUMLAUT e\else\protect \highumlaut e\fi
\else\if\string#1i\ifmmode\@MATHUMLAUT i\else\protect\highumlaut\i \fi
\else\if\string#1E\ifmmode\@MATHUMLAUT E\else\protect\highumlaut E\fi
\else\if\string#1I\ifmmode\@MATHUMLAUT I\else\protect\highumlaut I\fi
\else\if\string#1s\ifmmode\@SS\else\ss\fi{}%
\else\if\string#1-\allowhyphens\-\allowhyphens
\else\if\string#1\string"\hskip\z@skip
\else\if\string#1|\discretionary{-}{}{\kern.03em}%
\else\if\string#1c\allowhyphens\discretionary{k-}{}{c}\allowhyphens
\else\if\string#1f\allowhyphens\discretionary{ff-}{}{f}\allowhyphens
\else\if\string#1l\allowhyphens\discretionary{ll-}{}{l}\allowhyphens
\else\if\string#1m\allowhyphens\discretionary{mm-}{}{m}\allowhyphens
\else\if\string#1n\allowhyphens\discretionary{nn-}{}{n}\allowhyphens
\else\if\string#1p\allowhyphens\discretionary{pp-}{}{p}\allowhyphens
\else\if\string#1t\allowhyphens\discretionary{tt-}{}{t}\allowhyphens
\else\if\string#1<\flqq{}%
\else\if\string#1>\frqq{}%
\else             \dq #1%
\fi\fi\fi\fi\fi\fi\fi\fi\fi\fi\fi\fi\fi\fi\fi\fi\fi\fi\fi\fi\fi\fi\fi\fi\fi}
}
\gdef\dateaustrian{\def\today{\number\day.~\ifcase\month\or
  J\"anner\or Februar\or M\"arz\or April\or Mai\or Juni\or
  Juli\or August\or September\or Oktober\or November\or Dezember\fi
  \space\number\year}}
\gdef\dategerman{\def\today{\number\day.~\ifcase\month\or
  Januar\or Februar\or M\"arz\or April\or Mai\or Juni\or
  Juli\or August\or September\or Oktober\or November\or Dezember\fi
  \space\number\year}}
\gdef\dateUSenglish{\def\today{\ifcase\month\or
 January\or February\or March\or April\or May\or June\or
 July\or August\or September\or October\or November\or December\fi
 \space\number\day, \number\year}}
\gdef\dateenglish{\def\today{\ifcase\day\or
 1st\or 2nd\or 3rd\or 4th\or 5th\or
 6th\or 7th\or 8th\or 9th\or 10th\or
 11th\or 12th\or 13th\or 14th\or 15th\or
 16th\or 17th\or 18th\or 19th\or 20th\or
 21st\or 22nd\or 23rd\or 24th\or 25th\or
 26th\or 27th\or 28th\or 29th\or 30th\or
 31st\fi
 ~\ifcase\month\or
 January\or February\or March\or April\or May\or June\or
 July\or August\or September\or October\or November\or December\fi
 \space \number\year}}
\gdef\datefrench{\def\today{\ifnum\day=1\relax 1\/$^{\rm er}$\else
  \number\day\fi \space\ifcase\month\or
  janvier\or f\'evrier\or mars\or avril\or mai\or juin\or
  juillet\or ao\^ut\or septembre\or octobre\or novembre\or d\'ecembre\fi
  \space\number\year}}
%
%
\gdef\captionsgerman{%
\def\refname{Literatur}%
\def\abstractname{Zusammenfassung}%
\def\bibname{Literaturverzeichnis}%
\def\chaptername{Kapitel}%
\def\appendixname{Anhang}%
\def\contentsname{Inhaltsverzeichnis}%
\def\listfigurename{Abbildungsverzeichnis}%
\def\listtablename{Tabellenverzeichnis}%
\def\indexname{Index}%
\def\figurename{Abbildung}%
\def\tablename{Tabelle}%
\def\partname{Teil}}
\gdef\captionsenglish{%
\def\refname{References}%
\def\abstractname{Abstract}%
\def\bibname{Bibliography}%
\def\chaptername{Chapter}%
\def\appendixname{Appendix}%
\def\contentsname{Contents}%
\def\listfigurename{List of Figures}%
\def\listtablename{List of Tables}%
\def\indexname{Index}%
\def\figurename{Figure}%
\def\tablename{Table}%
\def\partname{Part}}
\gdef\captionsfrench{%
\def\refname{R\'ef\'erences}%
\def\abstractname{R\'esum\'e}%
\def\bibname{Bibliographie}%
\def\chaptername{Chapitre}%
\def\appendixname{Appendice}%
\def\contentsname{Table des mati\`eres}%
\def\listfigurename{Liste des figures}%
\def\listtablename{Liste des tables}%
\def\indexname{Index}%
\def\figurename{Figure}%
\def\tablename{Table}%
\def\partname{Partie}}%
\newcount\language 
\newcount\USenglish  \global\USenglish=0
\newcount\german     \global\german=1
\newcount\austrian   \global\austrian=2
\newcount\french     \global\french=3
\newcount\english    \global\english=4
\gdef\setlanguage#1{\language #1\relax
  \expandafter\ifcase #1\relax
  \dateUSenglish  \captionsenglish   \or
  \dategerman     \captionsgerman    \or
  \dateaustrian   \captionsgerman    \or
  \datefrench     \captionsfrench    \or
  \dateenglish    \captionsenglish   \fi}
\gdef\originalTeX{\mdqoff \umlauthigh
  \let\ss\original@ss \let\3\original@three
  \setlanguage{\USenglish}}
\gdef\germanTeX{\mdqon \umlautlow \let\ss\newss \let\3\ss
  \let\dospecials\german@dospecials
  \setlanguage{\german}}
} 
%
%
\germanTeX
%
%
\output={\if N\header\headline={\hfill}\fi
\plainoutput\global\let\header=Y}
\magnification\magstep1
\tolerance = 500
\hsize=14.4true cm
\vsize=22.5true cm
\parindent=6true mm\overfullrule=2pt
\newcount\kapnum \kapnum=0
\newcount\parnum \parnum=0
\newcount\procnum \procnum=0
\newcount\nicknum \nicknum=1
\font\ninett=cmtt9

\font\ninebf=cmbx9

\font\sixbf=cmbx6
\font\ninesl=cmsl9

\font\nineit=cmti9

\font\ninerm=cmr9

\font\sixrm=cmr6
\font\ninei=cmmi9
\font\eighti=cmmi8
\font\sixi=cmmi6
\skewchar\ninei='177 \skewchar\eighti='177 \skewchar\sixi='177
\font\ninesy=cmsy9
\font\eightsy=cmsy8
\font\sixsy=cmsy6
\skewchar\ninesy='60 \skewchar\eightsy='60 \skewchar\sixsy='60
\font\titelfont=cmr10 scaled 1440
\font\paragratit=cmbx10 scaled 1200

\font\name=cmcsc10
\font\emph=cmbxti10

\font\tenmsbm=msbm10
\font\sevenmsbm=msbm7
%

%
\font\got=eufm10

\font\teneufm=eufm10
\font\seveneufm=eufm7
\font\fiveeufm=eufm5
\newfam\eufmfam
\textfont\eufmfam=\teneufm
\scriptfont\eufmfam=\seveneufm
\scriptscriptfont\eufmfam=\fiveeufm

\font\tenmsam=msam10
\font\sevenmsam=msam7
\font\fivemsam=msam5
\newfam\msamfam
\textfont\msamfam=\tenmsam
\scriptfont\msamfam=\sevenmsam
\scriptscriptfont\msamfam=\fivemsam
\font\tenmsbm=msbm10
\font\sevenmsbm=msbm7
\font\fivemsbm=msbm5
\newfam\msbmfam
\textfont\msbmfam=\tenmsbm
\scriptfont\msbmfam=\sevenmsbm
\scriptscriptfont\msbmfam=\fivemsbm
\def\Bbb#1{{\fam\msbmfam\relax#1}}
\def\cz{{\kern0.4pt\Bbb C\kern0.7pt}
}
\def\ez{{\kern0.4pt\Bbb E\kern0.7pt}
}
\def\fz{{\kern0.4pt\Bbb F\kern0.3pt}}
\def\gz{{\kern0.4pt\Bbb Z\kern0.7pt}}
\def\hz{{\kern0.4pt\Bbb H\kern0.7pt}
}
\def\kz{{\kern0.4pt\Bbb K\kern0.7pt}
}
\def\nz{{\kern0.4pt\Bbb N\kern0.7pt}
}
\def\oz{{\kern0.4pt\Bbb O\kern0.7pt}
}
\def\rz{{\kern0.4pt\Bbb R\kern0.7pt}
}
\def\sz{{\kern0.4pt\Bbb S\kern0.7pt}
}
\def\pz{{\kern0.4pt\Bbb P\kern0.7pt}
}
\def\qz{{\kern0.4pt\Bbb Q\kern0.7pt}
}
\newskip\ttglue
\def\ninepoint{\def\rm{\fam0\ninerm}%
  \textfont0=\ninerm \scriptfont0=\sixrm \scriptscriptfont0=\fiverm
  \textfont1=\ninei \scriptfont1=\sixi \scriptscriptfont1=\fivei
  \textfont2=\ninesy \scriptfont2=\sixsy \scriptscriptfont2=\fivesy
  \textfont3=\tenex \scriptfont3=\tenex \scriptscriptfont3=\tenex
  \def\it{\fam\itfam\nineit}%
  \textfont\itfam=\nineit
  \def\sl{\fam\slfam\ninesl}%
  \textfont\slfam=\ninesl
  \def\bf{\fam\bffam\ninebf}%
  \textfont\bffam=\ninebf \scriptfont\bffam=\sixbf
   \scriptscriptfont\bffam=\fivebf
  \def\tt{\fam\ttfam\ninett}%
  \textfont\ttfam=\ninett
  \tt \ttglue=.5em plus.25em minus.15em
  \normalbaselineskip=11pt
  \font\name=cmcsc9
  \let\sc=\sevenrm
  \let\big=\ninebig
  \setbox\strutbox=\hbox{\vrule height8pt depth3pt width0pt}%
  \normalbaselines\rm
  \def\sl{\it}}

\headline={\ifodd\pageno\rightheadline\else\leftheadline\fi}
\def\rightheadline{\ninepoint Paragraphen"uberschrift\hfill\folio}
\def\leftheadline{\ninepoint\folio\hfill Chapter"uberschrift}
\let\header=Y
\def\titel#1{\need 9cm \vskip 2truecm
\parnum=0\global\advance \kapnum by 1
{\baselineskip=16pt\lineskip=16pt\rightskip0pt
plus4em\spaceskip.3333em\xspaceskip.5em\pretolerance=10000\noindent
\titelfont Chapter \uppercase\expandafter{\romannumeral\kapnum}.
#1\vskip2true cm}\def\leftheadline{\ninepoint
\folio\hfill Chapter \uppercase\expandafter{\romannumeral\kapnum}.
#1}\let\header=N
}
\def\Titel#1{\need 9cm \vskip 2truecm
\global\advance \kapnum by 1
{\baselineskip=16pt\lineskip=16pt\rightskip0pt
plus4em\spaceskip.3333em\xspaceskip.5em\pretolerance=10000\noindent
\titelfont\uppercase\expandafter{\romannumeral\kapnum}.
#1\vskip2true cm}\def\leftheadline{\ninepoint
\folio\hfill\uppercase\expandafter{\romannumeral\kapnum}.
#1}\let\header=N
}
\def\need#1cm {\par\dimen0=\pagetotal\ifdim\dimen0<\vsize
\global\advance\dimen0by#1 true cm
\ifdim\dimen0>\vsize\vfil\eject\noindent\fi\fi}
\def\neupara#1{\par\penalty-2000
\procnum=0\global\advance\parnum by 1
\vskip1cm\noindent{\paragratit \the\parnum. #1}%
\def\rightheadline{\ninepoint\S\the\parnum.\ #1\hfill \folio}%
\vskip 8mm\noindent}
\def\Proclaim #1 #2\finishproclaim {\bigbreak\noindent
{\bf#1\unskip{}. }{\it#2}\medbreak\noindent}
%
\gdef\proclaim #1 #2 #3\finishproclaim {\bigbreak\noindent%
\global\advance\procnum by 1
{%
{\relax\ifodd \nicknum
\hbox to 0pt{\vrule depth 0pt height0pt width\hsize
   \quad \ninett#3\hss}\else {}\fi}%
\bf\the\parnum.\the\procnum\ #1\unskip{}. }
{\it#2}
\immediate\write\num{\string\def
 \expandafter\string\csname#3\endcsname
 {\the\parnum.\the\procnum}}
\medbreak\noindent}
\newcount\stunde \newcount\minute \newcount\hilfsvar
\def\uhrzeit{
    \stunde=\the\time \divide \stunde by 60
    \minute=\the\time
    \hilfsvar=\stunde \multiply \hilfsvar by 60
    \advance \minute by -\hilfsvar
    \ifnum\the\stunde<10
    \ifnum\the\minute<10
    0\the\stunde:0\the\minute~Uhr
    \else
    0\the\stunde:\the\minute~Uhr
    \fi
    \else
    \ifnum\the\minute<10
    \the\stunde:0\the\minute~Uhr
    \else
    \the\stunde:\the\minute~Uhr
    \fi
    \fi
    }

 \def\calH{{\cal H}}

\def\calM{{\cal M}} \def\calN{{\cal N}}

 \def\calX{{\cal X}}
\def\calY{{\cal Y}} \def\calZ{{\cal Z}}

\def\gotp{\hbox{\got p}} 
\def\goto{\hbox{\got o}}

\def\dim{\mathop{\rm dim}\nolimits}

\def\GL{\mathop{\rm GL}\nolimits}

\def\id{\mathop{\rm id}\nolimits}

\def\kernel{\mathop{\rm kernel}\nolimits}

\def\mod{\mathop{\rm mod}\nolimits}
\def\O{{\rm O}}
\def\U{{\rm U}}

\def\SO{\mathop{\rm SO}\nolimits}

\def\Sp{\mathop{\rm Sp}\nolimits}

\def\boxit#1{
  \vbox{\hrule\hbox{\vrule\kern6pt
  \vbox{\kern8pt#1\kern8pt}\kern6pt\vrule}\hrule}}
\def\Boxit#1{
  \vbox{\hrule\hbox{\vrule\kern2pt
  \vbox{\kern2pt#1\kern2pt}\kern2pt\vrule}\hrule}}

\def\smallni{\smallskip\noindent }
\def\medni{\medskip\noindent }
\def\bigni{\bigskip\noindent }
\def\Isom{\mathop{\;{\buildrel \sim\over\longrightarrow }\;}}
\def\lo{\longrightarrow}

\def\loma{\longmapsto}
\def\betr#1{\vert#1\vert}
\def\spitz#1{\langle#1\rangle}
\def\pii{\pi {\rm i}}

\def\square{\hbox{\hbox to 0pt{$\sqcup$\hss}\hbox{$\sqcap$}}}
\def\qed{\ifmmode\square\else{\unskip\nobreak\hfil
\penalty50\hskip3em\null\nobreak\hfil\square
\parfillskip=0pt\finalhyphendemerits=0\endgraf}\fi}
\def\pn{\the\parnum.\the\procnum}
\def\downmapsto{{\buildrel
        {\vbox{\hbox{\hskip.2pt$\scriptstyle-$}}}
        \over{\raise7pt\vbox{\vskip-4pt\hbox{$\textstyle\downarrow$}}}}}

\nopagenumbers

\immediate\newwrite\num
\nicknum=0  
\let\header=N
\def\transpose#1{#1'}%
\def\Jac{{\rm Jac}}\def\tr{{\rm tr}}\def\triv{{\rm triv}}\def\sym{{\rm sym}}
\def\wt{{\rm wt}}
\def\Rho{\varrho_{\hbox{\sevenrm Jac}}}
\input vector-herm.num
\immediate\openout\num=vector-herm.num
\immediate\newwrite\num\immediate\openout\num=vector-herm.num
\def\RAND#1{\vskip0pt\hbox to 0mm{\hss\vtop to 0pt{%
  \raggedright\ninepoint\parindent=0pt%
  \baselineskip=1pt\hsize=2cm #1\vss}}\noindent}
\noindent
\centerline{\titelfont Vector valued hermitian and quaternionic modular forms}%
\def\leftheadline{\ninepoint\folio\hfill
Vector valued hermitian and quaternionic modular forms}%
\def\rightheadline{\ninepoint Introduction\hfill \folio}%
\headline={\ifodd\pageno\rightheadline\else\leftheadline\fi}
\vskip 1.5cm
\leftline{\it \hbox to 6cm{Eberhard Freitag\hss}
Riccardo Salvati
Manni  }
  \leftline {\it  \hbox to 6cm{Mathematisches Institut\hss}
Dipartimento di Matematica, }
\leftline {\it  \hbox to 6cm{Im Neuenheimer Feld 288\hss}
Piazzale Aldo Moro, 2}
\leftline {\it  \hbox to 6cm{D69120 Heidelberg\hss}
 I-00185 Roma, Italy. }
\leftline {\tt \hbox to 6cm{freitag@mathi.uni-heidelberg.de\hss}
salvati@mat.uniroma1.it}
\vskip1cm
\centerline{\paragratit \rm  2014}%
\vskip5mm\noindent%
\let\header=N%
\def\imag{{\rm i}}%
\def\St{{\rm St}}%
\vskip1.5cm\noindent
{\paragratit Introduction}%
\bigni
Extending the method of the paper [FS3] we prove three structure theorems
for vector valued modular forms, where two correspond to 4-dimensional cases
(two hermitian modular groups) and one to a 6-dimensional case (a quaternionic
modular group). In the literature there have been treated several cases of dimension
$\le 3$. We refer to [FS3] to some comments about what has been done in this area.
\smallskip
Vector valued modular forms belong to automorphy factors
which are related to rational
representations of the complexification of the maximal
compact group of the underlying Lie group. In the case of the symplectic group
$\Sp(n,\rz)$ 
this is $\GL(n,\cz)$. In the hermitian case (Lie group $\U(n,n)$)
it is $\GL(n,\cz)\times\GL(n,\cz)$
and in the quaternary case (Lie group $\Sp(n,\hz)$ or, in the notation
of Helgason, $\SO^*(4n)$) it is $\GL(2n,\cz)$.
We will describe these general facts briefly in this paper. More details can be
found in the Heidelberg diploma theses [Hey] (hermitian case)  and [SH]
(quaternionic case).
\smallskip
The basic automorphy
factor is the Jacobian. In the hermitian case it corresponds to the
external tensor product $\St\otimes\St$  of two standard representations of $\GL(n,\cz)$.
In the quaternionic case it corresponds to the 
representation 
$\Rho$ of $\GL(2n,\cz)$ with highest weight $(1,1,0,\dots,0)$.
\smallskip
We then restrict to $n=2$.
First we treat two hermitian modular groups,
one belonging to the field of Eisenstein numbers, the other to the field of 
Gaussian numbers.
In the Eisenstein case we take the congruence group of level
$\sqrt{-3}$ and in the Gauss case the congruence group of level $1+\imag$.
In both cases the groups have to be extended by the external involution $\tau(Z)=Z'$.
Their rings of modular forms have been determined in [Ma] and [FS2].
In both cases the ring of modular form is a polynomial ring in 5 generators of equal
weight. This implies that the corresponding modular varieties are projective spaces $P^4$.
So the two groups are very distinguished.
The main result in these cases states
that in both cases a certain module of  vector valued modular forms with respect
to the representations $\det^r\St\otimes\St$ is generated by 
Rankin--Cohen brackets $\{f,g\}$ among the generators of the rings of modular forms
(Theorems \MTa, \MTb).
\smallskip
In the quaternionic case we
treat a very particular modular group of degree two which belongs
to the ring of Hurwitz integers, more precisely to an extension of index 6
of the principal congruence subgroup of level $\gotp=(1+\imag_1)$. This group
has been introduced in [FH]. 
Its ring of modular forms has been determined in [FS1].
It is a weighted polynomial ring
where the generators have weights $3,1,1,1,1,1,1$.
The direct sum of the spaces of
vector valued modular forms with respect to $\det^r\Rho$ is a
module over this ring whose structure will be determined (Theorem \MTc).
Again generators are given by Rankin--Cohen brackets.
\neupara{The hermitian symplectic group}%
The hermitian symplectic group is the unitary group $\U(n,n)$.
We use the
matrix
$$J=\pmatrix{0&-E\cr E&0}\qquad (E\ \hbox{denotes the $n\times n$-unit matrix}).$$
The matrix $\imag J$ is hermitian of signature $(n,n)$. We
denote its unitary group by $\U(n,n)$. So this is the subgroup of $\GL(2n,\cz)$
defined by $\transpose{\bar M} JM=J$.
\smallskip
We denote by
$$\calH_n:=\{Z\in\cz^{n\times n};\quad \imag(\bar Z'-Z)>0\ \hbox{(positive definit)}\}$$
the hermitian half plane.
This is an open subset of the $\cz$-vector space $\cz^{n\times n}$. The group
$\U(n,n)$ acts on $\calH_n$ through the usual formula
$$MZ=(AZ+B)(CZ+D)^{-1},\quad M=\pmatrix{A&B\cr C&D}.$$
There is an extra automorphism $Z\mapsto Z'$ and one has 
$M(Z')'=\bar M(Z)$.
\smallskip
The following description of vector valued automorphy factors is due to Heyen [Hey].
The maps
$$\eqalign{
\U(n,n)\times\calH_n\lo\GL(n,\cz),\quad &(M,Z)\loma CZ+D,\cr
&(M,Z)\loma \bar CZ'+\bar D,\cr}$$
are both automorphy factors. Since they have equal rights we introduce
the group $\GL(n,\cz)\times\GL(n,\cz)$ (which is the complexification of the maximal
compact subgroup $\U(n)\times\U(n)$ of $\U(n,n)$).
If
$$\varrho:\GL(n,\cz)\times\GL(n,\cz)\lo\GL(\calZ)$$
is a rational representation on a finite dimensional complex vector space $\cz$,
we can consider the automorphic factor
$$\varrho(CZ+D,\bar C Z'+\bar D).$$
Usually we restrict to irreducible representations. They are tensor products
of two irreducible representations of $\GL(n,\cz)$. An irreducible rational representation
of $\GL(n,\cz)$ is called {\it reduced\/} if it is polynomial and if it does not vanish
along the determinant surface $\det(A)=0$. We can write every irreducible
rational representation
$\varrho$ of $\GL(n,\cz)\times\GL(n,\cz)$ in the form
$$\varrho(A,B)=\det A^{r_1}\det B^{r_2}\varrho_1(A)\otimes
\varrho_2(B)$$
with reduced representations $\varrho_i$.
\proclaim
{Lemma}
{Let be $M\in\U(n,n)$ and $Z\in\calH_n$. Then
$$\det(CZ+D)=\det M\det(\bar CZ'+\bar D).$$
}
DetAB%
\finishproclaim
The proof can be given by means of the usual generators of the group $\U(n,n)$ [Kr].
Details can be found in [Hey].\qed
Taking into account this Lemma, we are led to Heyen's definition
of automorphy factors
$$\det M^\epsilon\det(CZ+D)^r\varrho_1(CZ+D)\otimes\varrho_2(\bar C Z'+\bar D)$$
where $\epsilon$ and $r$ are integers and
where $\varrho_i$ are two reduced representations of $\GL(n,\cz)$.
These are automorphy factors for the whole Lie group $\U(n,n)$.
\neupara{Vector valued hermitian modular forms}%
The hermitian modular group with respect to an imaginary quadratic field $F$ is
$$\Gamma_F=\{M\in\U(n,n);\quad M\ \hbox{integral in}\ F\}.$$
In the following  $\Gamma\subset\Gamma_F$ denotes a subgroup of finite index.
We introduce vector valued modular forms for automorphy factors as defined in the
previous section. Since the discrete group $\Gamma$ may have more characters than the
power of the determinant, we should now 
consider instead of an integer $\epsilon\in\gz$ an arbitrary character 
$\chi:\Gamma\to\cz^*$. Let now
$r\in\gz$ be a weight and consider two reduced representations $\varrho_1,\varrho_2$
of $\GL(n,\cz)$. Then we can consider the representation $\varrho:=\varrho_1\otimes\varrho_2$
of $\GL(n,\cz)\times\GL(n,\cz)$.  We can consider the automorphy
factor
$$J(M,Z)=\chi(M)\det(CZ+D)^r\varrho(CZ+D,\bar CZ'+\bar D).$$
We denote by
$\calZ$ the representation space of $\varrho$.
A holomorphic modular form is a holomorphic function $f:\calH_n\to\calZ$
with the transformation property 
$$f(MZ)=J(M,Z)f(Z)\quad\hbox{for}\quad M\in\Gamma.$$ 
In the case $n=1$ the usual regularity
condition at the cusps has to be added. We denote this space by
$$[\Gamma,r,\chi,\varrho].$$
More generally, one can consider meromorphic solutions $f$ of the functional equation.
As meromorphicity condition at the cusps we demand that there exists a scalar valued
holomorphic modular form $g$ such that $gf$ is holomorphic. We denote the space of
meromorphic forms by
$$\{\Gamma,r,\chi,\varrho\}.$$
The field of modular functions is
$$K(\Gamma)=\{\Gamma,0,{{\rm triv}},{{\rm triv}}\}.$$
We want to study the following example. Let $\St=\id$ be the standard representation
of $\GL(n,\cz)$. We want to take $\varrho_1=\varrho_2=\St$. The representation
$\varrho=\St\otimes\St$ can
be realized on the space $\calZ=\cz^{n\times n}$ of $n\times n$-matrices by
$$\varrho(A,B)(W)=AW\transpose B.$$
Hence modular forms in this case are functions $f:\calH_n\to\cz^{n\times n}$ with the transformation
property
$$f(MZ)=\chi(M)\det(CZ+D)^r(CZ+D)f(Z)(\bar CZ'+\bar D)'.$$
We want to compute the derivative of a substitution $Z\mapsto MZ$.
The Jacobian of a substitution $M\in\U(n,n)$ at a point $Z\in\calH_n$
can be considered
as a linear map
$$\Jac(M,Z):\cz^{n\times n}\lo\cz^{n\times n}.$$
The following lemma is due to Pfrommer [Pf]. It 
can be proved by means of the standard generators of $\U(n,n)$.
\proclaim
{Lemma}
{The Jacobian $\Jac(M,Z)$ of a substitution $M\in\U(n,n)$ at a point $Z\in\calH$
is given by the linear map
$$\Jac(M,Z)(W)=(\bar CZ'+\bar D)'^{-1}W(CZ+D)^{-1}.$$
Its determinant is
$$\det\Jac(M,Z)=\det M^n\det(CZ+D)^{-2n}.$$
}
JacU%
\finishproclaim
We can consider the differentials
$dz_{ik}$ on $\calH_n$ and collect them to a matrix $dZ=(dz_{ik})$. From Lemma
\JacU\ we obtain the formula
$$dZ\vert M=(\bar CZ'+\bar D)'^{-1}dZ(CZ+D)^{-1}.$$
We want to consider meromorphic differentials on $\calH_n$. They can be written
in the form $\tr(fZ)$ where $f$ is an $n\times n$-matrix of meromorphic functions
on $\calH_n$. This differential is invariant under $\Gamma$ if and only if
$$f(MZ)=(CZ+D)f(Z)(\bar CZ'+\bar D)',$$
a formula which can be found already in [Pf].
This is the transformation law of  a form in $\{\Gamma,\chi,r,\varrho\}$
with $\chi(M)=1$, $r=0$ and $\varrho=\St\otimes\St$.
\neupara{Eisenstein numbers}%
Let
$$F=\qz[\omega],\quad \omega:={-1+\sqrt{-3}\over 2}$$
be the field of Eisenstein numbers.
Its ring of integers is $\gz[\omega]$.
In the paper [FS] the hermitian modular group of degree two
$$\Gamma_F=\U(2,2)\cap\GL(4,\goto_F)$$
and its congruence subgroup of level $\sqrt{-3}$
$$\Gamma[\sqrt{-3}]=\Gamma_F[\sqrt{-3}]:=\kernel(\Gamma_F\lo\GL(4,\gz[\omega]/\sqrt{-3}))$$
have been considered.
\smallskip
For integral $r$ we denote by $[\Gamma[\sqrt{-3}],r]$ the space of all scalar valued
holomorphic modular forms of transformation type
$$f(MZ)=\det M^r\det(CZ+D)^rf(Z)\qquad(M\in\Gamma[\sqrt{-3}]).$$
So in our vector valued notation this is
$$[\Gamma[\sqrt{-3}],r]=[\Gamma[\sqrt{-3}],r,{{\rm det}}^r,\triv].$$
We are interested in the subspace of all symmetric forms
$$[\Gamma[\sqrt{-3}],r]^\sym:=\{f\in[\Gamma[\sqrt{-3}],r];\ f(Z)=f(Z')\}.$$
We collect these spaces to the ring
$$A(\Gamma[\sqrt{-3}]):=\bigoplus_r[\Gamma[\sqrt{-3}],r]^\sym.$$
In [FS] it has been proved that
this ring is a polynomial ring in five theta constants
which have been introduced by Dern and Krieg in [DK]:
$$\Theta_p(Z):=\sum_{g\in\goto_F^2}e^{2\pi\hbox{\sevenrm i}
  \overline{(g+p)}'Z(g+p)},$$
where $\sqrt{-3}p$ runs through the five
$$\pmatrix{0\cr 0},\quad \pmatrix{1\cr 0},\quad
 \pmatrix{0\cr 1},\quad \pmatrix{1\cr 1},\quad
 \pmatrix{1\cr -1}.$$
These five series $\Theta_1,\dots,\Theta_5$
are algebraically independent modular forms of weight one. From [FS] we know the following
result.
\proclaim
{Theorem} 
{One has
$$A(\Gamma[\sqrt{-3}])=\cz[\Theta_1,\dots,\Theta_5].$$
Hence $A(\Gamma[\sqrt{-3}])$
is a polynomial ring in five variables.}
TheoFS%
\finishproclaim
We also consider the space of symmetric vector valued modular forms
$$\calM(r)=[\Gamma[\sqrt{-3}],r,{{\rm det}}^r,\St\otimes\St]^{\hbox{\sevenrm sym}}.$$
Recall that this are holomorphic functions $f:\calH_n\to\cz^{2\times 2}$
with the transformation property
$$f(MZ)=\det M^r\det(CZ+D)^r(CZ+D)f(Z)(\bar CZ'+\bar D)'\ \hbox{and}\ f(Z)=f(Z').$$
We collect them to
$$\calM=\bigoplus_r\calM(r).$$
This is a graded module over $A(\Gamma[\sqrt{-3}])$. 
\smallskip
Pfrommer [Pf] constructed examples of vector valued hermitian modular forms
by means of Rankin--Cohen brackets.
Following his construction
we obtain elements from
$\calM$ as follows. let $f,g\in[\Gamma[\sqrt{-3}],1]$. Then
$$\{f,g\}=g^2d(f/g)$$
can be considered as element of $\calM(2)$. 
In particular, we are interested in the
forms $\{\Theta_i,\Theta_j\}$.
\proclaim
{Lemma (Dern--Krieg)}
{If we consider the $4$ elements $\{\Theta_1,\Theta_i\}$, $1<i\le 5$, as $4\times 4$-matrix,
then its determinant equals up to constant factor $\Theta_1^3\phi_9$ where
$\phi_9$ is the skew symmetric modular form of weight $9$ with respect to the full
modular group which has been introduced by Dern and Krieg [DK].}
DetSym%
\finishproclaim
This lemma is just a reformulation of Corollary 6 in [DK].\qed
\proclaim
{Theorem}
{$$\calM=\sum_{1\le i<j\le 5}A(\Gamma[\sqrt{-3}])\{\Theta_i,\Theta_j\}.$$
Defining relations are
$$\Theta_k\{\Theta_i,\Theta_j\}=\Theta_j\{\Theta_i,\Theta_k\}+\Theta_i\{\Theta_k,\Theta_j\}.$$
}
MTa%
\finishproclaim
{\it Proof.\/} The space 
$$\{\Gamma[\sqrt{-3}],r,r,\St\otimes\St\}^{\hbox{\sevenrm sym}}$$
is a vector space over the field of symmetric modular functions of dimension $\le 4$
(which is is the rank of of the representation $\St\otimes\St$). We get a basis if we 
multiply $\{\Theta_1,\Theta_i\}$, $i>1$, by $\Theta_i^{r-1}$. Hence an arbitrary 
$T\in \{\Gamma[\sqrt{-3}],{{\rm det}}^r,r,\St\otimes\St\}^{\hbox{\sevenrm sym}}$ can be written in the form
$$T=\sum_{i=2}^5g_i\{\Theta_1,\Theta_i\}$$
where $g_i$ are scalar valued meromorphic modular forms. From Lemma \DetSym\ follows that
$\Theta_1^3\phi_9 g_i$ is holomorphic. Since this is a skew modular form, Corollary 7
of [DK] can be applied to show that that $\Theta_1^3g_i$ is holomorphic. In other words
$$\calM\subset {1\over\Theta_1^3}\sum_{j=2}^5A(\Gamma[\sqrt{-3}])\{\Theta_1,\Theta_j\}.$$
Since we can interchange the variables we get
$$\calM\subset\bigcap_i{1\over\Theta_i^3}\sum_jA(\Gamma[\sqrt{-3}])\{\Theta_i,\Theta_j\}.$$
Theorem \MTa\ is now an easy consequence.\qed
\neupara{Gauss numbers}%
Now we consider the Gauss number field $K=\qz[\imag]$, its hermitian modular group
$$\Gamma_K=\U(2,2)\cap\GL(4,\gz[i])$$
and its congruence group of level $1+\imag$
$$\Gamma[1+\imag]=\Gamma_K[1+\imag]:=\kernel(\Gamma_K\lo\GL(4,\gz[\imag]/(1+\imag))).$$
For {\it even\/} $r$ we denote by $[\Gamma[1+\imag],r]$ the space of all scalar valued
holomorphic modular forms of transformation type
$$f(MZ)=\det M^{r/2}\det(CZ+D)^rf(Z)\qquad(M\in\Gamma[1+\imag]).$$
So in our vector valued notation this is
$$[\Gamma[1+\imag],r]=[\Gamma[1+\imag],r,{{\rm det}}^{r/2},\triv]\qquad(r\equiv 0\mod 2).$$
We are interested in the subspace of all symmetric forms
$$[\Gamma[1+\imag],r]^\sym:=\{f\in[\Gamma[1+\imag],r];\ f(Z)=f(Z')\}.$$
We collect these spaces to the ring
$$A(\Gamma[1+\imag]):=\bigoplus_{r\equiv 0\mod 2}[\Gamma[1+\imag],r]^\sym.$$
Basic elements of these are  the ten theta series
$$\Theta[m]=\sum_{g\in\gz[\imag]^2}e^{\pii(\overline{(g+a/2)}'Z(g+a/2)+b'(g+a/2))},\quad
m=\pmatrix{a\cr b}.$$
Here $a,b$ are two columns  of the special form
$$a=(1+\imag)\alpha,\quad b=(1+\imag)\beta,\qquad
\alpha\in\{0,1\}^2,\ ,\beta\in\{0,1\}^2,\ \alpha'\beta\equiv 0\mod 2.$$
These theta functions have been introduced in [Fr]. The following result is due to Matsumoto
[Ma] and has been reproved by Hermann in [He].
\proclaim
{Theorem (Matsumoto)}
{We have 
$$A(\Gamma[1+\imag])=\cz[\dots\Theta[m]^2\dots].$$
The ten theta squares generate a five dimensional space. Hence $A(\Gamma[1+\imag])$
is a polynomial ring in five variables.
}
TheoMa%
\finishproclaim
For sake of completeness we mention that the following five ${\alpha\choose\beta}$ 
give generators:
$$\pmatrix{1&0&1&0&1\cr1&0&0&0&1\cr0&1&0&0&1\cr0&1&0&0&1\cr}$$
We denote the five corresponding thetas by $\Theta(1),\dots\Theta(5)$.
\smallskip
We consider (for even $r$) the space of symmetric vector valued modular forms
$$\calN(r)=[\Gamma[1+\imag],r,{{\rm det}}^{r/2},\St\otimes\St]^{\hbox{\sevenrm sym}}.$$
Recall that these are holomorphic functions $f:\calH_n\to\cz^{2\times 2}$
with the transformation property
$$f(MZ)=\det M^{r/2}\det(CZ+D)^r(CZ+D)f(Z)(\bar CZ'+\bar D)'\ \hbox{and}\ f(Z)=f(Z').$$
We collect them to
$$\calN=\bigoplus_{r\equiv0\mod2}\calN(r).$$
This is a graded module over $A(\Gamma[1+\imag])$. The forms
$$\{\Theta(i)^2,\Theta(j)^2\}=\Theta(j)^4d(\Theta(i)^2/\Theta(j)^2)$$
can be considered as element of $\calN(4)$.
\smallskip
We need two scalar valued modular forms $\phi_4$ and $\phi_{10}$ which do not belong
to the ring $A(\Gamma(1+\imag)$. The form $\phi_{10}$ has been introduced in[Fr]
as the product of the ten thetas $\Theta[m]$. It is a symmetric modular form of weight
$10$ with respect to the full modular group $\Gamma_K$ and the trivial multiplier system.
The form $\phi_4$ is a skew symmetric modular form of weight 4 with respect to the
full modular group and with trivial multiplier system too. It has been introduced
in [Ma] and occurs also in [He] and [DK]. 
\smallskip
The zero set of $\phi_4$ is the 
$\Gamma(1+\imag)$-orbit of the set defined by $Z=Z'$. The zero orders  are one.  Hence
every skew symmetric modular form with respect to $\Gamma(1+\imag)$ is divisible by $\phi_4$.
\smallskip
The zero sets of the single $\Theta[m]$ are also known. 
Since the full modular group 
permutes the 10 one-dimensional spaces generated by them transitively, it is sufficient
to treat one case. In the case $m'=(1,1,1,1)$ the function $\Theta[m]$ vanishes
along the fixed point set of the transformation
$$Z\loma Z'\left[\matrix{1&0\cr 0&\imag}\right]$$
in first order and every zero is equivalent mod $\Gamma[1+\imag]$ to a point of this
set.
\proclaim
{Lemma}
{If we consider the $4$ elements $\{\Theta(1)^2,\Theta(i)^2\}$, $1<i\le 5$, as $4\times 4$-matrix
then its determinant equals up to constant factor $\Theta_1^6\phi_4\phi_{10}$.}
DetSymG%
\finishproclaim
{\it Proof.\/} We denote the determinant by
$$D=\det(\{\Theta(1)^2,\Theta(2)^2\},\dots,\{\Theta(1)^2,\Theta(5)^2\}).$$
One can check that $D$ is skew symmetric and hence divisible by $\phi_4$.
Expanding the determinant one shows easily that $\Theta(1)^{-7}D$ is holomorphic.
The full modular group acts on the space generated by the five $\Theta(i)^2$.
From this it follows  that
$\Theta(1)^{-6}D$
is a modular form with respect to the full modular modular group. 
It is divisible by $\Theta(1)$ and hence by all $\Theta(i)$. So 
$D/(\Theta(1)^{-6}\phi_4\phi_{10})$ is a holomorphic modular form of weight $0$ and hence
constant.\qed
\proclaim
{Theorem}
{$$\calN=\sum_{1\le i<j\le 5}A(\Gamma[1+\imag])\{\Theta(i)^2,\Theta(j)^2\}.$$
Defining relations are
$$\Theta(k)^2\{\Theta(i)^2,\Theta(j)^2\}=\Theta(j)^2\{\Theta(i)^2,\Theta(k)^2\}+
\Theta(i)^2\{\Theta(k)^2,\Theta(j)^2\}.$$
}
MTb%
\finishproclaim
{\it Proof.\/}
Similar to the Eisenstein case
an arbitrary 
$$T\in \{\Gamma[1+\imag],{{\rm det}}^{r/2},r,\St\otimes\St\}^{\hbox{\sevenrm sym}}$$ 
can be written in the form
$$T=\sum_{i=2}^5g_i\{\Theta_1^2,\Theta_i^2\}$$
where $g_i$ are scalar valued meromorphic modular forms. 
We assume that $T$ is holomorphic. Then
from Lemma \DetSymG\ follows that
$\Theta(1)^6\phi_4\phi_{10} g_i$ is holomorphic. Since it is skew symmetric, it is divisible
by $\phi_4$ and we obtain that
$\Theta(1)^6\phi_{10} g_i$ is holomorphic. We want to show that this form is divisible
by $\phi_{10}$. It is sufficient to show that it is divisible by each $\Theta[m]$ and,
using the action of the full modular group, it is sufficient to restrict to
$m'=(1,1,1,1)$. The forms $\Theta(1)^2$ and $g_i$ are invariant under
$$Z\loma Z'\left[\matrix{1&0\cr 0&\imag}\right]$$
but $\phi_{10}$ changes it sign. Hence we can divide by $\Theta[m]$ and we obtain that
$\Theta(1)^6 g_i$ is holomorphic.
\smallskip
In other words
$$\calM\subset {1\over\Theta(1)^6}\sum_{j=2}^5A(\Gamma[1+\imag])\{\Theta(1)^2,\Theta(j)^2\}.$$
Since we can interchange the variables we get
$$\calM\subset\bigcap_i{1\over\Theta(i)^6}\sum_jA(\Gamma[1+\imag])\{\Theta(i)^2,\Theta(j)^2\}.$$
Theorem \MTb\ is now an easy consequence.\qed
\neupara{The quaternionic symplectic group}%
Let $\hz=\rz+\rz\imag_1+\rz\imag_2+\rz\imag_3$ be the field of quaternions.
The defining relations are
$$\eqalign{&\imag_1^2=\imag_2^2=\imag_3^2=-1,\cr
&\imag_1\imag_2=-\imag_2\imag_1=\imag_3,\ \imag_2\imag_3=-\imag_3\imag_2=\imag_1,\
    \imag_3\imag_1=-\imag_1\imag_3=\imag_2.\cr}$$
The conjugate $\bar x$ of a quaternion $x=x_0
+x_1\imag_1+x_2\imag_2+x_3\imag_3$ is
$$\bar x:=x_0-x_1\imag_1-x_2\imag_2-x_3\imag_3.$$
We use the embedding
$$\eqalign{
\hz\lo\cz^{2\times 2},\quad x&\loma\check x,\cr
x_0+x_1\imag_1+x_2\imag_2+x_3\imag_3&\loma
\pmatrix{x_0+\imag x_1&x_2+\imag x_3\cr -x_2+\imag x_3&x_0-\imag x_1}.\cr}$$
It induces an isomorphism of algebras
$$\hz\otimes_\rz\cz\Isom \cz^{2\times 2}.$$
We extend this to an isomorphism of matrix algebras
$$(\hz\otimes_\rz\cz)^{n\times n}\Isom \cz^{2n\times 2n}$$
by applying the check operator componentwise. In particular, we obtain an
embedding
$$\GL(n,\hz)\lo \GL(2n,\cz),\quad A\loma \check A.$$
This is a complexification in the sense that the Lie algebra of $\GL(2n,\cz)$ arises
as complexification of the Lie algebra of $\GL(n,\hz)$.
\smallskip
The quaternionic symplectic group $\Sp(n,\hz)$ consists
of all $2n\times 2n$-matrices $M$ with entries in $\hz$ such that
$$\bar M'JM=J,\qquad J=\pmatrix{0&-E\cr E&0}\quad (E\ \hbox{unit matrix}).$$
(This is the group $\SO^*(4n)$ in the notation of Helgason.)
Let $\calX_n\subset\hz^{n\times n}$ be the set of quaternionic
hermitian $n\times n$ matrices
and $\calY_n$ the cone of all positive definit ones. The quaternionic half plane is
$$\calH_n=\calX_n+\imag\calY_n\subset(\hz\otimes_\rz\cz)^n.$$
This is an open subset of the complex vector space
$$\calZ_n=\calX_n\otimes_\rz\cz=\calX_n+\imag\calX_n.$$
Its dimension is $2+2n(n-1)$.
\smallskip
The group $\Sp(n,\hz)$ acts on $\calH_n$ through the usual formula
$$Z\loma  MZ=(A+BZ)(CZ+D)^{-1},\quad M=\pmatrix{A&B\cr C&D}.$$
The standard generators of $\Sp(n,\hz)$ and their actions are
\smallni
a) $\pmatrix{E&H\cr 0&E}(Z)=Z+H\qquad (\bar H'=H\in\calX_n)$.\hfill\break
b) $\pmatrix{U&0\cr 0&\bar U'^{-1}}(Z)=UZ\bar U'\qquad (U\in\GL(n,\hz))$.\hfill\break
c) $J(Z)=-Z^{-1}$.
\medni
In the case b) we can apply the transformation to all $Z\in\calZ_n$. This gives a
representation
$$\varrho:=\GL(n,\hz)\lo\GL(\calZ_n),\quad U\loma (W\mapsto UW\bar U').$$
\proclaim
{Remark}
{There exists a unique rational (holomorphic)
representation $\Rho:=\GL(2n,\cz)\lo\GL(\calZ_n)$
with the property
$$\Rho(\check U)=\varrho(U)\quad\hbox{for}\quad U\in\GL(n,\hz).$$}
RatExt%
\finishproclaim
{\it Proof.\/} The derived representation of $\varrho$ extends to a (complex)
representation of the Lie algebra of $\GL(2n,\cz)$. Each representation of the
Lie algebra of $\GL(m,\cz)$ on a finite dimensional complex vector space is induced
by a rational representation of $\GL(m,\cz)$.\qed
\smallskip
We can consider the Jacobian
$${{\rm Jac}}:\Sp(n,\hz)\times \calH_n\lo\GL(\calZ_n).$$
\proclaim
{Lemma}
{We have
$${{\rm Jac}}(M,Z)=\Rho(\check C\check Z+\check D)^{-1}.$$
{\bf Corollary.}
$$\det{{\rm Jac}}(M,Z)=\det(\check C\check Z+\check D)^{-3}.$$
}
JacRho%
\finishproclaim
{\it Proof.\/} Since both sides are automorphy factors, it is sufficient to prove
this for generators. For the translations  $Z\mapsto Z+H$ and
unimodular transformations $Z\mapsto \bar U'ZU$ the statement is trivial. 
Hence it is sufficient to treat the case $M=J$. The Jacobian can be computed
easily.
\smallni
{\it The Jacobian of the transformation $Z\mapsto -Z^{-1}$  is the
linear map $W\mapsto Z^{-1}WZ^{-1}$.}
\smallni
Since the  claimed formula in Lemma \JacRho\ is an identity between rational 
functions in $Z\in\calZ_n$ it is sufficient to prove it 
for a {\it real\/} invertible matrix $Z\in\calX_n$.
For real $Z$ we have 
$\Rho(\hat Z)^{-1}=\varrho(Z)^{-1}$.
This is the linear map $W\loma Z^{-1}WZ^{-1}$.
\qed
\smallskip
From Lemma \JacRho\ we get a different description of the representation $\Rho$.
The stabilizer $\Sp(n,\hz)_{\imag E}$ of the point $\imag E$.
It consists of matrices of the form
$$\pmatrix{A&-B\cr B&A}.$$
We attach to such a matrix $\check A+\imag\check B$.
This defines an ismorphism
$$\Sp(n,\hz)_{\imag E}\lo \U(2n).$$
Hence the map
$$\Sp(n,\hz)_{\imag E}\lo \GL(2n,\cz).$$
is a complexification map. The assignment $M\mapsto \Jac(M,\imag E)$ is a representation
of $\Sp(n,\hz)_{\imag E}$. Its complexification is $\Rho$.
\proclaim
{Lemma}
{The representation
$$\Rho:\GL(2n,\cz)\lo\GL(\calZ_n)\qquad (\dim_\cz\calZ_n=2+2n(n-1))$$
is irreducible. Its highest weight is $(1,1,0,\dots,0)$.}
HigW%
\finishproclaim
{\it Proof.\/} Recall that a highest weight vector is a vector in $\calZ_n$
that is invariant under all unipotent upper-triangular matrices.
If $A\in\GL(n,\hz)$ is a unipotent upper-triangular matrix, then $\check A$
has the same property. Hence highest weight vectors have the property
$\varrho(A)a=a$. It follows that the matrix $E_{11}\in\calZ_n$ whose $(1,1)$-entry
is 1 and all other entries are 0 is up to a constant factor the only highest
weight vector. Hence the representation is irreducible.
\smallskip
The statement about the highest weight says that a diagonal matrix $D\in\GL(2n,\cz)$ acts
on the highest weight vector by multiplication by $d_{11}d_{22}$. Since
the representation is holomorphic, it is sufficient to proof this for
diagonal matrices of special shape. We can assume that the diagonal elements
are of the form $d_1,\bar d_1,\dots,d_n,\bar d_n$.
Then the matrix has to act by multiplication by $\betr {d_1}^2$.
Write $d_i=\alpha_i+\imag\beta_i$. Consider in $\GL(n,\hz)$ the diagonal matrix $D_1$
with diagonal entries $\alpha_i+\imag_1\beta_i$. Then $\check D_1=D$.
Now the stated formula can be verified.\qed
\neupara{A ring of quaternionic modular forms of degree two}%
We denote by $\Omega_n$ the
image of $\Sp(n,\hz)$ in the group of biholomorphic self maps
of $\calH_n$. Since two symplectic matrices define the same transformation
if and only if they agree up to the sign, we have
$$\Omega_n=\Sp(n,\hz)/\{\pm E\}.$$
The map
$$j:\Sp(n,\hz)\times\calH_2\lo\cz^*,\quad (M,Z)\loma\det(\check C\check Z+\check D),$$
is a scalar valued factor of automorphy. Since $j(M,Z)=j(-M,Z)$
it factors through an automorphy factor of $\Omega_n$,
$$j:\Omega_n\times\calH_n\lo\cz^*.$$
In the case $n>2$ the group $\Omega_n$ is the full group of biholomorphic self maps
of $\calH_n$. But in the case $n=2$ the full group is extension of index two
$$\hat\Omega_2=\Omega\cup\tau\Omega,\quad \tau(Z)=Z'.$$
The automorphy factor $j$ can be extended to an automorphy factor
$$j:\hat\Omega_2\times\calH_2\lo\calH_2,\quad j(\tau,Z)=1.$$
We denote by $\goto$ the ring of Hurwitz integers. It contains
$$\goto_0:=\gz+\gz\imag_1+\gz\imag_2+\gz\imag_3$$
as subgroup of index two which is extended
by the element
$$\omega:={1+\imag_1+\imag_2+\imag_3\over2}.$$
We consider the two sided ideal $\gotp$ which can be generated by any
element of norm 2 as for example $1+\imag_1$.
\smallskip
In the paper [FH] a certain subgroup
$$\Gamma(\goto)[\gotp]\subset\hat\Omega_2$$ has been defined. It is related to the
principal congruence subgroup
$$\Sp(2,\goto)[\gotp]=\kernel(\Sp(2,\goto)\lo\Sp(2,\goto/\gotp)).$$
We need a certain extension of this group. For this we need the following type of
transformations. Let $\sigma:\hz\to \hz$ be an $\rz$-linear transformation which is
orthogonal in the sense that it preserves the standard scalae product $\bar xx$.
We extend $\sigma$ to hermitian matrices in $\calX_2$ by
$$\sigma\pmatrix{x_0&x_1\cr *&x_2}=
\pmatrix{x_0&\sigma(x_1)\cr *&x_2}$$
Since positive definite matrices are mapped to positive definit ones, we can extend
$\sigma$ to a biholomorphic transformation of $\calH_2$. Hence the orthogonal group
$\O(2,\rz)$ acts on $\calH_2$. This defines a certain subgroup of $\hat\Omega_2$.
The group of diagonal matrices in $\Sp(2,\hz)$ produces the subgroup
$\SO(2,\rz)$ and the extension by $\tau$ produces $\O(2,\rz)$.
\proclaim
{Definition}
{The group $\O'(4,\gz)$
consists of all transformations
$$(x_1,x_2,x_3,x_4)\loma(\pm
    x_{\sigma(1)},\pm x_{\sigma(2)},
            \pm x_{\sigma(3)},\pm x_{\sigma(4)}),$$
where $\sigma$ is an arbitrary permutation and
where the number of minus-signs is even.}
Oplus%
\finishproclaim
This is a subgroup of
index $2$ in $\O(4,\gz)$.
\proclaim
{Definition}
{The group
$\Gamma(\goto)[\gotp]$ is defined as
$$
    \Gamma(\goto)[\gotp]=
    \spitz{\Sp(2,\goto)[\gotp]/\{\pm E\},\O'(4,\gz)}.$$}
LeveL%
\finishproclaim
This group has been introduced in [FH]. Using the exceptional isogeny between
$\O(2,6)$ and $\Sp(2,\hz)$ in this paper also a very natural orthogonal description
of this group has been given.
\proclaim
{Lemma}
{The group
$\Gamma(\goto)[\gotp]$
contains $\Sp(2,\goto)[\gotp]/\{\pm E\}$ as normal subgroup of index $6$.
The quotient group is isomorphic
to $S_3$. It can be generated by the images of
$$x\loma\omega x\bar\omega,\qquad x\loma\alpha\bar x\bar\alpha\quad\hbox{where}
\quad\alpha={1+\imag_1\over\sqrt2}.$$
}
NormaL%
\finishproclaim
We consider modular forms for this group. The space $[\Gamma(\goto)[\gotp],r]$ for an
integral $r$ consists of all holomorphic functions $f:\calH_2\to\cz$ with the transformation
property
$$f(\gamma(Z))=j(\gamma,Z)^rf(Z)\quad\hbox{for all}\quad\gamma\in [\Gamma(\goto)[\gotp].$$
We collect them to the graded algebra
$$A(\Gamma(\goto)[\gotp])=\bigoplus_{r\in\gz}[\Gamma(\goto)[\gotp],r].$$
The structure of this algebra has been investigated in [FH] and, correcting an error
in [FH], it has finally been determined in [FS]. In [FS] the orthogonal language has been
used. The translation between the orthogonal and the symplectic picture can be found
in detail in [FH].
\proclaim
{Theorem}
{The algebra $A(\Gamma(\goto)[\gotp])$ is a polynomial ring of $7$ variables,
$G,\vartheta_1,\dots,\vartheta_6$ where
$G$ is of weight $3$ and the other six of weight $1$.}
AlStru%
\finishproclaim
We mention that the weights of the symplectic picture are doubled
in the orthogonal picture. Hence in the paper [FS] the weights are $6,2,\dots,2$.
\neupara{Vector valued quaternionic modular forms of degree 2}%
In section one we studied the vector valued automorphy factor Jac. It factors
through $\Omega_n$ and, in the case $n=2$, extends to $\hat\Omega_2$.
$$\Jac:\hat\Omega_2\times\calH_2\loma\GL(\calZ_2).$$
We denote by $\calM(r)$ the space of all holomorphic functions
$$f:\calH_2\lo\calZ_2$$
with the transformation property
$$f(\gamma(Z))=j(\gamma,Z)^r\Jac(\gamma,Z)f(Z)\quad\hbox{for all}\quad
\gamma\in\Gamma(\goto)[\gotp].$$
The direct sum
$$\calM=\bigoplus_{r\in\gz}\calM(r)$$
is a graded module over $A=A(\Gamma(\goto)[\gotp])$.
\smallskip
For two non-zero homogenous elements of positive degree $f,g\in A$ we define
$$\{f,g\}:=\wt(g)g df-\wt(f)fdg.$$
Here $\wt(f)$ denotes the weight of $f$.
Another way to write this is
$$\{f,g\}={g^{\wt(f)+1}\over f^{\wt(g)-1}}\; d\Bigl(
{f^{\wt(g)}\over g^{\wt(f)}}\Bigr).$$
This shows
$$\{f,g\}\in\calM(\wt(f)+\wt(g)).$$
This is a skew-symmetric $\cz$-bilinear pairing and it satisfies
the following rule
$$\wt(h)h\{f,g\}=\wt(g)g\{f,h\}+\wt(f)f\{h,g\}.$$
We identify the vector space $\calZ_2$ with $\cz^6$ using the coordinates
$z_0$, $z_2$, $z_{10}$, $z_{11}$, $z_{12}$, $z_{13}$ where
$$Z=\pmatrix{z_0&z_1\cr *&z_2},\quad z_1=z_{10}+\imag_1z_{11}+
\imag_2z_{12}+\imag_3z_{13}.$$
Then we can consider the brackets $\{f,g\}$ as columns with 6 entries.
\proclaim
{Lemma}
{$$\det(\{G,\vartheta_1\},\{\vartheta_2,\vartheta_1\},\dots,\{\vartheta_6,\vartheta_1\})=
\vartheta_1^{5}D.$$
Here $D$ is modular form of weight $12$ with respect to a non-trivial character
(more precisely: $D^2\in[\Gamma(\goto)[\gotp],24]$).
The zero divisor of $D$
is the $\Gamma(\goto)[\gotp]$-orbit of the set $z_{10}=z_{11}$. The vanishing order
is one.
}
LemJac%
\finishproclaim
{\it Proof.\/}
We denote by
$\nabla g$ the  column of the partial derivatives of a holomorphic function $g$ on
$\calH_2$ with respect to the variables $z_0,\dots,z_{13}$
In [FS] it has been proved that the determinant
$$D=\det\pmatrix{\vartheta_1&\dots&\vartheta_6&3G\cr \nabla
\vartheta_1&\dots&\nabla \vartheta_6&\nabla G}$$
is a modular form of weight 12 with respect to a nontrivial  character.
Now
$$\vartheta_1^6 D =\det\pmatrix{\vartheta_1&\vartheta_1
\vartheta_2&\dots&\vartheta_1\vartheta_6&\vartheta_1 3G\cr 
\nabla\vartheta_1&\vartheta_1 \nabla\vartheta_2
&\dots&\vartheta_1\nabla \vartheta_6&\vartheta_1\nabla G}$$
Let us  multiply the first  column , by $\vartheta_2$  and   
subtract  it to  the second column,  and so on.
Now
$$\vartheta_1^6 D =\det\pmatrix{\vartheta_1&0&\dots&0&0\cr 
\nabla\vartheta_1& \{\vartheta_1,\vartheta_2\} &\dots& \{\vartheta_1,\vartheta_6\} 
& \{\vartheta_1,G\}}$$
Hence
$$\vartheta_1^5 D =\det\pmatrix {\{\vartheta_1,\vartheta_2\} 
&\dots& \{\vartheta_1,\vartheta_6\} & \{\vartheta_1,G\}}.$$
In [FH] (Lemma 2.1 in connection with Proposition 8.4) has been proved
that the zero divisor
is the $\Gamma(\goto)[\gotp]$-orbit of the set $z_{10}=z_{11}$. The vanishing order
is one. Lemma \LemJac\ is an immediate consequence.\qed
\proclaim
{Theorem}
{The $A(\Gamma(\goto)[\gotp])$-module $\calM$ is generated by the
brackets $\{f,g\}$ between the $7$ generators. Defining relations are
$$\{f,g\}=-\{g,f\}\quad\hbox{and}\quad \wt(h)h\{f,g\}=\wt(g)g\{f,h\}+\wt(f)f\{h,g\}$$
for any  $f,g,h$ of the $7$ generators
}
MTc%
\finishproclaim
{\it Proof.\/}
Since the determinant of the matrix $(\{G,\vartheta_1\},\{\vartheta_2,\vartheta_1\},\dots,\{\vartheta_6,\vartheta_1\})$
does not vanish, every $T\in\calM(r)$ can be written in the form
$$T=g_1\{G,\vartheta_1\}+g_2\{\vartheta_2,\vartheta_1\}+\cdots+g_6\{\vartheta_6,\vartheta_1\}$$
with meromorphic functions $g_i$. They transform as modular forms in the ring
$A(\Gamma(\goto)[\gotp])$.
From Lemma \LemJac\ we see that $g_i\vartheta_1^{5}D$ is holomorphic.
Since $D$ vanishes along $z_{10}=z_{11}$ in first order and since the transformation
$z_{10}\leftrightarrow z_{11}$ is in the group $\Gamma(\goto)[\gotp]$ the forms
$g_i\vartheta_1^{5}D$ are skew symmetric with respect to this transformation.
The forms $g_i$, $\vartheta_1$ are symmetric. Hence $g_i\vartheta_1^{5}D$ vanishes along
$z_{10}=z_{11}$ and must be divisible by $D$. So we see
$$T\in{1\over \vartheta_1^{5}}\bigl(A\{G,\vartheta_1\}+A\{\vartheta_2,\vartheta_1\}
+\cdots+A\{\vartheta_6,\vartheta_1\}\bigr).$$
One can take any $\vartheta_i$ instead of $\vartheta_1$ and obtains
$$\calM\subset\bigcap_{i=1}^6{1\over \vartheta_i^{5}}\Bigl( A\{G,\vartheta_i\}+
\sum_{j\ne i}A\{\vartheta_j,\vartheta_i\}\Bigr).$$
Theorem \MTc\ is an easy consequence.\qed
\vskip1.5cm\noindent
{\paragratit References}
\bigskip
\item{[DK]} Dern, T., Krieg, A.:
{\it Graded rings of Hermitian modular forms of degree 2,\/}
Manuscripta math. {\bf 110}, 251–-272 (2003)
\medskip
\item{[FH]} Freitag, E. Hermann, C.F.: {\it Some modular
varieties in low dimension,\/} Advances in Math. {\bf 152}, 203--287 (2000)
\medskip
\item{[Fr]} Freitag, E.: {\it Modulformen zweiten Grades zum rationalen und Gau"s'schen
Zahlk"orper,\/} Sitzungsber.\ Heidelb.\ Akad.\ d.\ Wiss. {\bf 1}, 1--49 (1967)
\medskip
\item{[FS1]} Freitag, E., Salvati Manni, R.:
{\it Hermitean modular forms and the Burkhardt quartic,\/} 
 Manuscripta math. {\bf 119}, 57--59 (2006) 
\medskip
\item{[FS2]} Freitag, E.,  Salvati Manni, R.:
{\it Some modular varieties of low dimension~II,\/}
Advances in Math. {\bf 214}, no.1, 132--145 (2007)
\medskip
\item{[FS3]} Freitag, E.,  Salvati Manni, R.:
{\it Basic vector valued
Siegel modular forms of genus two,\/}
preprint (2013)
\medskip
\item{[He]} Hermann, C.F.: {\it Some modular varieties related to $\pz^4$,\/}
Abelian Varieties, edited by Barth/Hulek/Lange, Walter de Gruyter Berlin-New York
(1995)
\medskip
\item{[Hey]} Heyen, D.: {\it Vektorwertige Hermite'sche Modulformen zweiten Grades zum
Gau"s'schen Zahlk"orper,\/}
Diplomarbeit am Mathematischen Institut der Universit"at Heidelberg (2008)
\medskip
\item{[Kr1]} Krieg,\ A.:
{\it Modular Forms on Half-Spaces of Quaternions, \/}
Lecture Notes in Math.~{\bf 1143}, Springer Verlag
Berlin-Heidelberg-New York Tokyo (1985)
\medskip
\item{[Kr2]} Krieg, A.: {\it The graded ring of quaternionic modular
forms of degree 2,\/} Math.\ Z. {\bf251}, 929--942 (2005)
\medskip
\item{[Ma]} Matsumoto, K.: {\it Thetafunctions on the bounded symmetric domain if type
$I_{2,2}$ and the period map of a 4-parameter family of $K_3$-surfaces,\/}
Math.\ Ann.\ {\bf 295}, 383--409 (1993)
\medskip
\item{[Pf]} Pfrommer,\ T.: {\it Konstruktion vektorwertiger Hermite'scher Modulformen
zweiten Grades,\/} 
Diplomarbeit am Mathematischen Institut der Universit"at Heidelberg (2010)
\medskip
\item{[SH]} Schulte--Hengesbach, F.:
{\it The six dimensional Hermitian Domain of Type IV.\/}
Diplomarbeit am Mathematischen Institut der Universit"at Heidelberg (2006)
\bye